\input amstex.tex
\loadmsam
\loadmsbm
\loadbold
\input amssym.tex
\input amstex

\input amssym.def
\input amssym.tex


\def\item#1{\vskip1.3pt\hang\textindent {\rm #1}}


\tolerance=300
\pretolerance=200
\hfuzz=1pt
\vfuzz=1pt


\hoffset=0.6in
\voffset=0.8in

\hsize=5.8 true in 


\vsize=8.5 true in
\parindent=25pt
\mathsurround=1pt
\parskip=1pt plus .25pt minus .25pt
\normallineskiplimit=.99pt

\countdef\revised=100
\mathchardef\emptyset="001F 
\chardef\ss="19
\def\3{\ss}
\def\anf{$\lower1.2ex\hbox{"}$}
\def\frac#1#2{{#1 \over #2}}
\def\>{>\!\!>}
\def\<{<\!\!<}

\def\ssarr{\hbox to 30pt{\rightarrowfill}}
\def\sarr{\hbox to 40pt{\rightarrowfill}}
\def\arr{\hbox to 60pt{\rightarrowfill}}
\def\larr{\hbox to 60pt{\leftarrowfill}}
\def\Arr{\hbox to 80pt{\rightarrowfill}}

{}



%
%


\def\Ind{\mathop{\rm Ind}\nolimits}






\def\0{{\bf 0}}
\def\1{{\bf 1}}

\def\a{{\frak a}}

\def\C{{\Bbb C}}

\def\R{{\Bbb R}}

\def\:{\colon}  
\def\.{{\cdot}}
\def\|{\Vert}
\def\bsk{\bigskip}

\def\giantskip{\vskip2\bigskipamount}
\def\gsk{\giantskip}

\def\msk{\medskip}

\def \res {\!\mid\!\!}

\def\bbr{\bigbreak}
\def\giantbreak{\par \ifdim\lastskip<2\bigskipamount \removelastskip
         \penalty-400 \giantskip\fi}

\def\nin{\noindent}
\def\cen{\centerline}
\def\pagebreak{\vskip 0pt plus 0.0001fil\break}
\def\linebreak{\break}

\def\hat{\widehat}

\def\phi{\varphi}
\def\epsilon{\varepsilon}
\def\eset{\emptyset}

\def\nin{\noindent}

\def\pder#1,#2,#3 { {\partial #1 \over \partial #2}(#3)}
\def\pde#1,#2 { {\partial #1 \over \partial #2}}


\def\subeq{\subseteq}

\def\tilde{\widetilde}

\def\up{{\uparrow}}

\font\eightrm=cmr8


\font\smc=cmcsc10
\font\bfone=cmbx10 scaled\magstep1 
\font\bftwo=cmbx10 scaled\magstep2 



\def\Lemma #1. {\bigbreak\vskip-\parskip\noindent{\bf Lemma #1.}\quad\it}

\def\Sublemma #1. {\bigbreak\vskip-\parskip\noindent{\bf Sublemma #1.}\quad\it}

\def\Proposition #1. {\bigbreak\vskip-\parskip\noindent{\bf Proposition #1.}
\quad\it}

\def\Corollary #1. {\bigbreak\vskip-\parskip\nin{\bf Corollary #1.}
\quad\it}

\def\Theorem #1. {\bigbreak\vskip-\parskip\noindent{\bf Theorem #1.}
\quad\it}

\def\Definition #1. {\rm\bigbreak\vskip-\parskip\noindent{\bf Definition #1.}
\quad}

\def\Remark #1. {\rm\bigbreak\vskip-\parskip\noindent{\bf Remark #1.}\quad}

\def\Example #1. {\rm\bigbreak\vskip-\parskip\noindent{\bf Example #1.}\quad}

\def\Problems #1. {\bigbreak\vskip-\parskip\noindent{\bf Problems #1.}\quad}
\def\Problem #1. {\bigbreak\vskip-\parskip\noindent{\bf Problems #1.}\quad}

\def\Conjecture #1. {\bigbreak\vskip-\parskip\noindent{\bf Conjecture #1.}\quad}

\def\Proof#1.{\rm\par\ifdim\lastskip<\bigskipamount\removelastskip\fi\smallskip
            \noindent {\bf Proof.}\quad}

\def\Axiom #1. {\bigbreak\vskip-\parskip\noindent{\bf Axiom #1.}\quad\it}

\def\Satz #1. {\bigbreak\vskip-\parskip\noindent{\bf Satz #1.}\quad\it}

\def\Korollar #1. {\bbr\vskip-\parskip\nin{\bf Korollar #1.} \quad\it}

\def\Bemerkung #1. {\rm\bigbreak\vskip-\parskip\noindent{\bf Bemerkung #1.}
\quad}

\def\Beispiel #1. {\rm\bigbreak\vskip-\parskip\noindent{\bf Beispiel #1.}\quad}
\def\Aufgabe #1. {\rm\bigbreak\vskip-\parskip\noindent{\bf Aufgabe #1.}\quad}

\def\Beweis#1. {\rm\par\ifdim\lastskip<\bigskipamount\removelastskip\fi
           \smallskip\noindent {\bf Beweis.}\quad}

\nopagenumbers

\def\date{\ifcase\month\or January\or February \or March\or April\or May
\or June\or July\or August\or September\or October\or November
\or December\fi\space\number\day, \number\year}

\def\title{Title ??}
\def\author{Author ??}

\def\thanks#1{\footnote*{\eightrm#1}}

\def\rightheadline{\hfil{\eightrm\title}\hfil\tenbf\folio}
\def\leftheadline{\tenbf\folio\hfil{\eightrm\author}\hfil}
\headline={\vbox{\line{\ifodd\pageno\rightheadline\else\leftheadline\fi}}}

\def\firstheadline{}
\def\firstfootline{\cen{\rm\folio}}

\def\seite #1 {\pageno #1
               \headline={\ifnum\pageno=#1 \firstheadline
               \else\ifodd\pageno\rightheadline\else\leftheadline\fi\fi}
               \footline={\ifnum\pageno=#1 \firstfootline\else{}\fi}}

\newdimen\dimenone
 \def\checkleftspace#1#2#3#4{
 \dimenone=\pagetotal
 \advance\dimenone by -\pageshrink   
 \ifdim\dimenone>\pagegoal          
   \else\dimenone=\pagetotal
        \advance\dimenone by \pagestretch
        \ifdim\dimenone<\pagegoal
          \dimenone=\pagetotal
          \advance\dimenone by#1         
          \setbox0=\vbox{#2\parskip=0pt                
                     \hyphenpenalty=10000
                     \rightskip=0pt plus 5em
                     \noindent#3 \vskip#4}    
        \advance\dimenone by\ht0
        \advance\dimenone by 3\baselineskip   
        \ifdim\dimenone>\pagegoal\vfill\eject\fi
          \else\eject\fi\fi}


\def\lsectionheadline #1 #2{\bigbreak\vskip-\lastskip
      \checkleftspace{1.1cm}{\bf}{#1}{\bigskipamount}
         \vbox{\vskip1.1cm}\cen{\bfone #1}\msk \cen{\bfone #2}\bsk}

\def\lchapterheadline #1 #2{\bigbreak\vskip-\lastskip\indent\vskip3cm
                       \cen{\bftwo #1} \msk \cen{\bftwo #2} \gsk}
\def\llsectionheadline #1 #2 #3{\bigbreak\vskip-\lastskip\indent\vskip1.8cm
\cen{\bfone #1} \msk \cen{\bfone #2} \msk \cen{\bfone #3} \nobreak\bsk\nobreak}


\newtoks\literat
\def\[#1 #2\par{\literat={#2\unskip.}%
\hbox{\vtop{\hsize=.15\hsize\nin [#1]\hfill}
\vtop{\hsize=.82\hsize\nin\the\literat}}\par
\vskip.3\baselineskip}

\mathchardef\emptyset="001F 
\def\address{Author: \tt$\backslash$def$\backslash$address$\{$??$\}$}


\pageno=1

\baselineskip=13pt plus 2pt
\documentstyle{amsppt}
\pageheight{45pc}
\pagewidth{33pc}


\def\up#1{\leavevmode \raise.16ex\hbox{#1}}





\chardef\ss="19

\def\3{\ss}

\def\bx{$\blacksquare$}
\def\Ind{\,\text{Ind}\,}

\def\lno{\leqalignno}
 
\def\ov{\overline}
\def\bs{\backslash}

\overfullrule=0pt
 
\def\Box #1 { \msk\par\nin 

\centerline{

\vbox{\offinterlineskip

\hrule

\hbox{\vrule\strut\hskip1ex\hfil{\smc#1}\hfill\hskip1ex}

\hrule}\vrule}\msk }

\def\bs{\backslash}

\topmatter
\title The Bochner measure and holomorphic extension of elementary spherical functions\endtitle
\author B. Kr\" otz (*), R. A. Kunze and R. J. Stanton (**)\endauthor
\leftheadtext{Kr\"otz, Kunze and Stanton}
\rightheadtext{The Bochner measure for spherical functions}
\address 
\noindent Max Planck Institute\hfil\break University of Georgia (ret.),\hfil\break
Ohio State University\hfil\break
\endaddress
\thanks (*) This research was supported in part by a Heisenberg grant of the DFG \endthanks
\thanks (**) This research was supported in part by NSF grant DMS 0301133 \endthanks
\endtopmatter
\vfill\eject

\document

\subhead  Introduction\endsubhead
\medskip

Let $G$ be a noncompact simple linear Lie group, and let $K$ be a maximal compact subgroup. Let 
$\hat G$ denote the set of equivalence classes of irreducible unitary representations of $G$ and let 
$\hat G_s$ denote those containing a trivial $K$-type. We recall several results from the work of 
Harish-Chandra. Any $(\pi,\Cal H_{\pi}) \in \hat G_s$ contains the trivial $K$-type with 
multiplicity one; let $v_0\in \Cal H_{\pi}$ be one of norm one. Denote by $\varphi_{\pi}$ the matrix 
coefficient
$$\varphi_{\pi}(g) = \langle \pi(g)v_0,v_0\rangle.$$
Then $\varphi_{\pi}$ is in $C^\omega(G//K), \varphi_{\pi}(e) =1, \text {and } \varphi_{\pi}$ is a 
positive definite function on $G$. One can give an explicit realization of such $(\pi,\Cal H_{\pi})$ 
by using induced representations, and in this way enlarge the parameter domain for these functions. 
So let $P = MAN$ be a minimal parabolic subgroup with $M\subset K$. Set $\pi_{\sigma, 
\lambda}=\Ind_{P}^G (\sigma\otimes\chi_{\lambda}\otimes\bold1)$ acting on ${\Cal 
H}_{\sigma,\lambda}$ where $(\sigma, W_\sigma) \in \hat M$, and $\chi_{\lambda}\in \text{Hom}(A,\Bbb 
C^x)\cong \frak a^\ast_{\Bbb C}$. For $\sigma =\bold1, \pi_\lambda:=\pi_{\bold1,\lambda}$ contains the 
trivial $K$ type with multiplicity one. While, for $\lambda \in \frak a^\ast_{\Bbb C}$, 
$\pi_\lambda$ need not in general be irreducible nor unitarizable, we do have $\hat 
G_s\hookrightarrow \{(\pi_\lambda, {\Cal H}_{\lambda})|\lambda \in \frak a^\ast_{\Bbb C}\}$ and, 
after incorporating the intertwining operator into the inner product, this gives unitary 
realizations. The explicit form of the inner product on ${\Cal H}_{\lambda}$ gives
$$\phi_\lambda(g) :=\phi_{\pi_\lambda}(g) = \int_K a(gk)^{\lambda-\rho}\ dk$$ 
which is clearly entire on $\frak a^\ast_{\Bbb C}$. If we set 
$$\varphi(x:\lambda) := \phi_{\lambda}(x), x\in G, \lambda \in \frak a^\ast_{\Bbb C}$$
then $\varphi$ is real analytic in ``$x$"and entire in ``$\lambda$".

In [KS-I] we constructed a universal domain in the complexification, $G_\Bbb C$, of $G$ and we 
showed that the matrix coefficients obtained from $K$-finite vectors of any irreducible Banach 
representation of $G$ have holomorphic extension to this domain. When one of the vectors is the 
trivial $K$-type this domain naturally maps to one in $G_{\Bbb C}/K_{\Bbb C}$. In [K{\'O}S] we 
showed that this domain, denoted $\Xi$, is a maximal $G$-invariant domain with this property, 
although not a maximal $K_\Bbb C$-invariant one. Consequently $\varphi(\cdot:\lambda)$ has 
holomorphic extension to $\Xi$. This domain is determined mainly by imposing two constraints, namely 
the holomorphic continuation of the $A$ component from the Iwasawa decomposition and the 
construction of non-unitary principal series in conjunction with the subrepresentation theorem. In 
this note we shall give an alternative justification of the holomorphic extension in the specific 
setting of spherical unitary representations and elementary spherical functions. The novelty of this 
analysis includes the definition of a {\it spectral} version of the Abel transform $f\to F_f$, and 
an integral formula which we make explicit for split rank one groups and from which we derive 
several consequences.

The starting point is the comment above that $\varphi_{\pi}$ is a positive definite function. The 
characterization of continuous positive definite functions on $\Bbb R$ by means of their Fourier
transforms and more importantly their relevance to such areas as spectral theory and probability was 
discovered by Bochner.  While there have been extensions of the characterization theorem to locally 
compact abelian groups and to some extent non-abelian groups, they have not played a central role in 
the harmonic analysis on these groups.  Yet a central theme in the development of harmonic analysis 
on semisimple Lie groups has been to examine the objects of harmonic analysis--matrix coefficients 
and characters--on certain abelian (or almost abelian) subgroups.  It seems natural then, following 
Bochner's lead, to try to characterize the abelian Fourier transforms of these objects and hopefully 
to find some uses of them in harmonic analysis.  This paper is the result of our first look at this 
problem.

The setting which seemed most amenable to this investigation is the non-unitary principal series of 
representations induced from a minimal parabolic subgroup $MAN$.  Here the abelian group $A$ enters 
in three well-known ways:  part of the starting data is a non-unitary character of $A$; the 
representation space can be realized on functions on a subgroup $V$ on which $A$ acts in a 
reasonable way; the growth properties of $K$-finite matrix coefficients are determined by their 
restrictions to $A$. 

In order to apply abelian harmonic analysis on $A$ in the most direct way to matrix coefficients of 
principal series representations, we shall use a minor modification of the ``non-compact picture" of 
non-unitary principal series which is built from an $A$-invariant measure on $V$.  The measure is 
constructed in $\S 3$ and the realization in $\S 4$. The Fourier transform, i.e., the Bochner 
measure, of matrix coefficients (not necessarily $K$-finite) of unitary principal series follows 
easily from this realization and is included in $\S 4$.  In $\S 5$ we 
consider the $\Bbb R$-rank 1 case and show that the Bochner measure can be computed for non-unitary 
principal series, $K$-finite coefficients provided $\vert\text{Re}~\lambda\vert<\vert\rho\vert.$  
For the zonal spherical functions, in $\S 6$ we compute the Bochner measure explicitly for split rank one 
groups, obtaining an integral 
representation which seems new even from a classical special function point of view.  In $\S 7$ we 
give four applications of the Bochner measure. In one we derive the Harish-Chandra expansion directly from 
the integral formula without any differential equations. Another is to define a {\it spectral} $F_f$ with 
properties in the spectral parameter analogous to those of $F_f$ in the group parameter.
\vfil
\eject

\subhead \S 1. Laplace transforms \endsubhead
\medskip
In classical abelian harmonic analysis there is a well known relationship, due to Paley and Wiener, 
between the holomorphic extension of functions to tubular neighborhoods and the exponential decay of 
their Fourier-Laplace transform. This will be examined in the context of spherical functions. We 
shall identify spherical functions as Laplace transforms on $\hat A\cong i\a^*\cong\a^*$. 
For this purpose we briefly review some standard facts on Laplace transforms. 
\msk Let $V$ be a finite dimensional real vector space and $V^*$ its dual. We denote by 
$V_\C=V+ iV$ the complexification of $V$. Let $\mu$ be a
positive Borel measure on $V^*$ and set 
$$D_\mu =\{x\in V\: \int_{V^*}  e^{-\alpha(x)}\ d\mu(\alpha)<\infty\}.$$
Shortly we shall assume that $\mu$ is {\it admissible}, i.e., $D_\mu\neq \eset$ (thus, in 
particular, that $\mu$ is a Radon measure). The {\it Laplace transform of} 
$\mu$ is then defined by 
$${\Cal L}\mu\: D_\mu\to \R_0^+, \ \ x\mapsto \int_{V^*} e^{-\alpha(x)}\ d\mu(\alpha).$$   
Note that ${\Cal L}\mu$ automatically extends to a holomorphic 
function on the tube domain $T_{D_\mu} =iV +\text{int }D_\mu\subeq V_\C$. We denote the extension also by 
${\Cal L}\mu$. We shall need the following results on Laplace transforms.

\Proposition 1.1. Let $\mu$ be an admissible positive Borel measure on the dual $V^*$ 
of a finite dimensional real vector space $V$. Then  
\item{(i)} the domain of definition $D_\mu$ of the Laplace transform ${\Cal L}\mu$ 
is a convex subset of $V$; 
\item{(ii)} suppose that $\mu$ is finite and that there is an open connected set $\Omega\subeq V$ 
with $0\in \Omega$ such that ${\Cal L}\mu$ has an extension to a
holomorphic function on the tube domain $T_\Omega=iV+\Omega$, then we 
have $\Omega\subeq  D_\mu$, i.e., on $T_\Omega$ the extension is represented by 
$$T_\Omega\to\C, \ \ z\mapsto \int_{V^*} e^{-\alpha(z)} \ d\mu(\alpha).$$ 
\Proof. (i) [Ne, Prop.\ V.4.3].
\par\nin (ii) [Ri, p.\ 311].\hfill\bx
 
\msk Let $A=\R^n$ be a simply connected abelian real Lie group. We identify the unitary dual $\hat 
A$ of $A$ with $(\R^n)^*$ by means of the isomorphism, 
$$(\R^n)^*\to \hat A, \ \ \alpha\mapsto (x\mapsto e^{i\alpha(x)}).$$ 
Let ${\Cal B}(\hat A)$ denote the Borel $\sigma$-algebra on $\hat A$ and for $(\pi, {\Cal H})$
a unitary representation of $A$ let $P({\Cal H})$ be the set of selfadjoint projections on ${\Cal 
H}$. Then there exists a spectral measure

$$E\:\ \ \Cal B (\hat A)\to P({\Cal H})$$
such that 
$$(\forall x\in A)\ \ \pi(x)=\int_{\hat A} e^{i\alpha(x)} \ dE(\alpha).$$

\msk We are now going to apply the above considerations to spherical functions
i.e. as matrix coefficients of the $K$-fixed vector of spherical principal series representations. 
Recall from [KS-I] the definition of the set $\Omega\subset\frak a$ (or see below) and from Prop. 
4.1 therein that $\phi_\lambda$ holomorphically extends to $G exp(i\Omega) K_{\Bbb C}$, and from 
Theorem 4.2 to $K_{\Bbb C}exp (i2\Omega) K_{\Bbb C}$.

\Lemma 1.2. Let $\lambda\in i\a^*$ and $\phi_\lambda$ be the associated spherical 
function on $G/K$. Then there exists a Radon probability measure $\mu$ on $\a^*$ 
such that 
$$(\forall a\in exp\, i2\Omega) \qquad \phi_\lambda(a^{-1})=\int_{\a^*} e^{i\alpha(\log a)} \ 
d\mu(\alpha).$$ 
\Proof. We consider the unitary representation 
$(\pi, {\Cal H}) =(\pi_\lambda\res_A, {\Cal H}_\lambda)$ of $A$. The probability measure $\mu$ in 
question is then given by  $\mu=E_{v_0, v_0}$ with $v_0=\bold 1_K$. 
In particular, we get that 
$$(\forall a\in A)\qquad \phi_\lambda(a^{-1})=\int_{\a^*} e^{i\alpha(\log a)}\ d\mu(\alpha).$$
Now the assertion of the Theorem follows from 
[KS-I Theorem 4.2] and Proposition 1.1(ii).\hfill\bx

\Remark  1.3. For $\lambda\in i\a^*$, as mentioned before, the spherical function $\phi_\lambda$ is 
a continuous positive definite function with $\phi_\lambda (e) =1$. In particular on the abelian 
group $A\cong \a$, $\phi_\lambda\res A$ has these properties. Hence the existence of a probability 
measure $\mu$ on $\a^*$ such that 

$$\phi_\lambda(a^{-1})=\int_{\a^*} e^{i\alpha(\log a)}\ d\mu(\alpha)$$

\noindent holds for all $a\in A$ is also a consequence of Bochner's Theorem. (In fact, Bochner used 
his result to give a proof of the spectral theorem.) However, the validity of the extension to 
$exp\, 2i\Omega$ of the Bochner integral representation is far from obvious. For example, the 
geometry of the tube reflects growth estimates of the Bochner measure and these are usually 
difficult to compute. All these remarks are illustrated in considerable detail in the case of rank 1 
groups in subsequent sections.

In view of Lemma 1.2 and Proposition 1.1(ii) every restricted spherical function
$\phi_\lambda\circ \roman{ exp_A}$, has a holomorphic extension  
to a tube domain

$$T_{\lambda, { \roman max}}=\a +i\omega_\lambda$$
which is maximal with respect to 
$\omega_\lambda\subeq \a$ open, convex, $0\in \omega_\lambda$ and $\omega_\lambda$ is ${\Cal 
W}_\a$-invariant. In [K-S I Th. 4.2] we show that $2\Omega \subset\omega_\lambda$ and in \S 5  for 
all real rank one groups $2\Omega =\omega_\lambda$ where $2\Omega=\{ X\in\a\:(\forall \alpha\in 
\Sigma)\  |\alpha(X)|<\pi\}.$

\subhead \S2.  The Technique \endsubhead
\medskip

Let $X$ be a locally compact space and let $A$ be a locally compact group with a continuous proper action 
$A \times X\to X$, denoted $(a,x)\to a\cdot x$.  Let $dx$ be a measure on $X$ and $da$ a left Haar
measure on $A$.  We shall make some assumptions to specialize the presentation for the simple
applications we have in mind.

First concerning the action, we shall assume that $A$ acts without fixed points so that each orbit 
is homeomorphic 
to $A$; we shall require that the orbit space $A\bs X$ is paracompact; and we want the measure $dx$ to 
be $A$ invariant. We will be 
concerned only with the case in which $A$ is abelian and isomorphic to a finite dimensional real 
vector 
space. Under such conditions, there is a unique measure 
$d\omega$ on $A\bs X$ and the usual ``double integration" formula for reasonable $f$ ([B] p. 44). 

$$\int\limits_Xf(x) dx=\int\limits_{A\bs X}\int\limits_Af(a\cdot x) da~d\omega.$$ 

Moreover, if we make $A\times A\bs X$ an $A$-space by translation in the $A$ coordinate, then $X$ and 
$A\times A\bs X$ are homeomorphic $A$-spaces ([B] p. 73).  If $f$ is a function on $X$, by abuse of 
notation, we will let $f$ also denote the corresponding function on $A\times A\bs X.$  The integration 
formula then becomes
$$\int\limits_X f(x)~dx=\int\limits_{A\times A\bs X}f(a, \omega) da\times d\omega$$

Since each orbit type is basically $A$ we can transfer $A$ harmonic analysis to each orbit.  Thus, 
if $\hat{A}$ denotes the unitary character group of $A$ with dual measure, for reasonable functions 
$f$ and $g$ on $X$ we have the familiar abelian results:

$$\lno{\text{if}~\chi_\lambda\in\hat{A},\, \hat{f}(\lambda, 
\omega)=\int\limits_Af(t,\omega)\chi_\lambda(t)^{-1}
da;&&(1)'\cr}$$
$\text{for}~a\in A~\text{with the usual transport of structure set}~L^\ast_af(t, \omega)=f(a^{-1}t,\omega),~\text{then}$
$$\lno{\widehat{L^\ast_af}(\lambda,\omega)=\chi_\lambda(a^{-1})\hat{f}(\lambda,\omega);&&(2)'\cr}$$

\noindent and the Parseval formula 

$$\lno{\int\limits_Af(t,\omega)\ov{g(t,\omega)}da=\int\limits_{\hat{A}}\hat{f}(\lambda,\omega)\ov
{\hat{g}(\lambda,\omega)}d\lambda.&&(3)'\cr}$$

For our purposes, we shall need to consider the usual action $L^\ast_a$ on functions twisted by a 
character.  This causes minor changes in $(1)', (2)'~\text{and}~ (3)'$.  So suppose $\chi_\upsilon$ 
is a continuous
character of $A$ (not necessarily unitary) and define

$$a \cdot f(t,\omega)=\chi_\upsilon(a)L^\ast_a f(t, \omega).$$
For the function $f_\upsilon$ defined by $f_\upsilon(t,\omega)=\chi_\upsilon(t^{-1})f(t,\omega)$ we 
have

$$\eqalignno{(a\cdot f)_\upsilon(t,\omega)&=\chi_\upsilon(t^{-1})\chi_\upsilon(a)L^\ast_a f(t, 
\omega)\cr
&=\chi_\upsilon(a^{-1}t)^{-1}L^\ast_a f(t,\omega)\cr
&=L^\ast_a f_\upsilon(t,\omega).\cr}$$
We shall call the contragredient to $\chi_\upsilon$ the character $\chi_{\upsilon'}$ defined by
$\chi_\upsilon\ov{\chi_{\upsilon'}}=1.$  For this $A$-action the modified Fourier transform defined 
by

$$\lno{\tilde{f}(\lambda,\omega):=\hat{f}_{\upsilon}(\lambda,\omega)&&(2.1)\cr}$$
has the property

$$\lno{\widetilde{a\cdot f}(\lambda, 
\omega)=\chi_\lambda(a^{-1})\tilde{f}(\lambda,\omega),&&(2.2)\cr}$$
and the Parseval formula
$$\lno{\int\limits_A f(a, \omega)\overline{g(a, \omega)}da =\int\limits_{\hat{A}}\hat{f}_{\upsilon}
(\lambda, \omega)\overline{\hat{g}_{\upsilon'}(\lambda, \omega)}d\lambda&&(2.3)\cr}$$

$$\qquad=\int\limits_{\hat{A}}\tilde{f}(\lambda, \omega)\overline{\tilde{g}(\lambda, 
\omega)}d\lambda$$

\noindent provided $f_\upsilon$ and $g_{\upsilon'}$ satisfy reasonable growth conditions.
\vfil
\eject

\subhead \S 3.  The invariant measure on V \endsubhead
\medskip

An interesting example of such an $A$-space (in fact our motivation) arises in semisimple Lie 
groups.  
We shall briefly recall some basic facts about these groups.  (See [H1]).

Let $G$ be a real semisimple Lie group with Lie algebra $\frak g$.  Fix a Cartan involution and let 
$\frak g=
\frak k \oplus \frak p$ be the corresponding decomposition.  Choose a maximal abelian subspace
$\frak a \subseteq \frak p.$  The adjoint action of $\frak a$ on $\frak g$ is diagonalizable; let 
$\Sigma$ 
denote the set of non-zero real eigenvalues and $\frak g_\alpha, \alpha \in \Sigma$, the corresponding 
eigenvectors.  
Fix a 
lexicographic order on $\Sigma$ and let $\Sigma^+$ be the positive elements.  If subalgebras of 
$\frak g$ 
are defined by $\frak n=\sum\limits_{\alpha\in\Sigma^+}\frak g_\alpha, \frak 
v=\sum\limits_{-\alpha\in\Sigma^+}\frak g_\alpha$ and $\frak m_1=$ zero eigenspace, then $\frak 
g=\frak n \oplus \frak m_1\oplus 
\frak v.$  Let $V, N, M_1, A$ denote the corresponding subgroups.  Then  $M_1=MA$ where $M\subset K$ 
and centralizes $A$.  Moreover, $V, N$ and $A$ are homeomorphic to their respective Lie algebras.  

Since $\frak v$ consists of $\frak a$-eigenspaces, the group $A$ acts on it and hence on $V$.  Denote 
this action on $V$ (conjugation) by $(a, v)\to a\cdot v$ or $\delta_a(v),$ and on functions by $\delta_a^\ast f(v) = f(\delta_{a^{-1}}(v)).$ There are two aspects of 
this 
action which we must examine before using $\S 2$:  (1) the obvious measure on $V$, Haar measure 
$dv,$ is 
not invariant but only relatively invariant; (2) the action of $A$ on $V$ is not free.  We shall 
remedy the first and define away the second.

The multiplier for the $A$ action on $dv$ is given by the character $\chi_{-2\rho}$, where 
$2\rho=\sum
\limits_{\alpha \in\Sigma^+}\alpha,$ i.e. for $f\in C_c(V)$

$$\delta_{a\ast}dv[f]:=dv[\delta_a^\ast f]=\chi_{-2\rho}(a) dv[f]=\chi_{2\rho}(a^{-1})dv[f].$$
We will use $\chi_{-\rho}$ to construct an invariant measure on $V$.

The eigenspace decomposition of $\frak g$ has an analog for $G$, namely the Bruhat big cell, such 
that except for a 
closed set of lower dimension $G=VMAN$ with the unique $A$-component denoted $a(g).$  Let $w\in G$ 
represent
the Weyl group element that sends $\Sigma^+$ to $-\Sigma^+.$  Then $w^2$ represents the identity,
$wNw^{-1}=V$ and for $v\neq e, w^{-1}v\in $ $VMAN$ ([Kn-S]).  We extend $\chi_{-\rho}$ to a function 
$\xi_{-\rho}$ defined a.e. on $G$ by 
$$\xi_{-\rho}(x):=\chi_{-\rho}(a(x)).$$

\proclaim {Lemma 3.1}  If $x$ is in $VMAN$ and $a$ is in $A$ then $ax$ and $xa$ are in $VMAN$.  
Also, 
$\xi_{-\rho}(ax)=\chi_{-\rho}(a)\xi_{-\rho}(x)$ and $\xi_{-\rho}(xa)=\xi_{-\rho}(x)\chi_{-\rho}(a).$
\endproclaim

\demo{Proof}  Obvious. \hfill\bx\enddemo
\medskip

Define $\chi_{w\rho}(a)=\chi_\rho(w^{-1}aw).$  Then $\chi_{w\rho}$ is a character and 
 
\proclaim{Lemma 3.2}  $\chi_{w\rho}$ is contragredient to $\chi_\rho.$  \endproclaim

\demo{Proof}  Since $\chi_\rho$ is $\Bbb R$-valued, this follows from $w\Sigma^+=\Sigma^-$.   \hfill\bx \enddemo

\proclaim{Proposition 3.3}  The measure $\xi_{-\rho}(w^{-1}v) dv$ on $V$ is $A$-invariant.  
\endproclaim

\demo{Proof}  Let $f\in C_c(V).$  First, notice that

$$\eqalignno{\delta^\ast_a(\xi_{-\rho}(w^{-1}\cdot)f)(v)&=\xi_{-\rho}(w^{-1}a^{-1}va) f(a^{-1}\cdot 
v)\cr
&=\chi_{-\rho}(w^{-1}a^{-1}w) \xi_{-\rho}(w^{-1}v)\chi_{-\rho}(a) \delta_a^\ast f(v)\cr
&=\chi_{-w\rho}(a^{-1})\chi_{-\rho}(a)\xi_{-\rho}(w^{-1}v)\delta_a^\ast f(v)\cr
&=\chi_{-2\rho}(a)\xi_{-\rho}(w^{-1}v)\delta_a^\ast f(v).\cr}$$
Here we have used Lemmas 3.1 and 3.2.  Now

$$\eqalignno{\delta_{a^\ast}(\xi_{-\rho}(w^{-1}\cdot) dv)[f]&=\xi_{-\rho}(w^{-1}\cdot) 
dv[\delta_a^\ast f]\cr
&=dv[\xi_{-\rho}(w^{-1}\cdot)\delta_a^\ast f]\cr
&=dv[\chi_{2\rho}(a)\delta_a^\ast(\xi_{-\rho}(w^{-1}\cdot) f)]\cr
&=\chi_{2\rho}(a)\delta_{a^\ast} dv[\xi_{-\rho}(w^{-1}\cdot)f]\cr
&=\chi_{2\rho}(a)\chi_{-2\rho}(a) dv [\xi_{-\rho}(w^{-1}\cdot) f]\cr
&=\xi_{-\rho}(w^{-1}) dv [f].\cr}$$\hfill\bx\enddemo

To obtain a space on which $A$ acts freely we shall discard a finite number of lower dimensional 
subspaces 
of $\frak v$ leaving an open dense subset $\frak v'$ which can be identified with a similar set 
$V'\subseteq 
V.$  For the order $\Sigma^+$ let $\alpha_1, \cdots \alpha_\ell$ be a basis of simple restricted 
roots.  
For each $i$ let $W_i=\sum\limits_{t\alpha_i\in\Sigma^+}\oplus\frak g_{t\alpha_i}$.  Each $W_i$ has 
a unique 
complementary subspace which is a sum of eigenspaces, say $W_i^c$.  Let 
$W^c=\mathop{\cup}\limits_{i=1}^\ell 
W_i^c$ and set $\frak v'=\frak v\bs W^c.$

\proclaim{Lemma 3.4}  $A$ acts freely on $\frak v'$.  \endproclaim

\demo{Proof}  If for each $i, P_i$ denotes projection along $W_i^c$ of $\frak v\to W_i$, then $\frak 
v'$ 
can be characterized as $\mathop{\cap}\limits_{i=1}^\ell \{ x\in \frak v\vert P_i x\neq 0\}.$  Now 
each 
$a \in A$ has all non-zero eigenvalues and commutes with each $P_i,$ therefore $\frak v'$ is 
$A$-invariant. 
Also, $\alpha_1, \cdots, \alpha_\ell$ form a basis for the dual space $\frak a'$, so for each $a \in 
A, a\neq 1,$ there is an $i$ with $\chi_{\alpha_i}(a)\neq 1.$  If $x\in \frak v'$ were fixed by $a$ then 
$P_ix=
P_i a\cdot x=a\cdot P_i x.$  But $P_i x\neq 0$ and on $W_i$ the element $a$ acts by 
$\chi_{\alpha_i}(a)$ or $\chi_{\alpha_i}
(a)^2$ both of which are not 1.  Thus the action is free.  \enddemo

We let $V'$ be the corresponding dense open set in $V$ and restrict $\xi_{-\rho}(w^{-1}v) dv$ to 
$V'$.  This 
serves as the $(X, dx)$ of $\S 2$.  

\remark{Remarks}  (1)  If dim $\frak a=1$ then $\frak v'=\frak v \bs \{ 0\}$ and so is maximal.  If 
dim $\frak a>1$ 
the maximal subset on which $A$ acts freely is not necessarily $\frak v'$ but can be constructed.  
As we have 
no use for it here and the construction is standard, we omit it.

(2)  $M$ also acts on $V$ by $(m, x)\to m\cdot x=mxm^{-1}.$  One can show that $M$ acts on $V'$ and 
the 
measure $\xi_{-\rho}(w^{-1}v) dv$ is $M$-invariant.  However, for the elementary spherical functions 
we shall not need this.

(3) It is worth noting that all that is done here can also be done for degenerate principal series 
induced off a maximal parabolic with $\frak n$ a prehomogeneous vector space of parabolic type. 
\endremark
\vfil
\eject

\subhead \S 4.  The Bochner measure - unitary induction \endsubhead
\medskip

We recall briefly the construction of the principal series representations of $G$ induced from the 
minimal parabolic subgroup $MAN$.  We shall use as Hilbert space for these representations, the 
``noncompact picture" consisting of vector valued square integrable functions on $V$.  Since the 
$A$-invariant
measure is absolutely continuous with respect to Haar measure on $V'$, a standard procedure provides 
an alternate realization of the principal series.

Let $(\sigma, W_\sigma)\in \hat M$, and let $\chi_\lambda$ be a 
character of $A$, not necessarily unitary.  In the ``induced picture" one uses

$$C^\infty(G; \sigma \otimes\chi_\lambda)=\{ f:  G\to W_\sigma\mid f~\text{is}~C^\infty~\text{and }
f(gman)=\sigma(m)^{-1}\chi_{\lambda+\rho}(a)^{-1}f(g)\}$$
and has $G$ act by left translation.  From the Iwasawa decomposition, $f$ is determined by 
its values on $K$ and one takes as norm the $L^2$ norm of $f$ restricted to $K$.  From the Bruhat 
big cell, $f$ is also determined by its values on $V$.  In this ``non-compact picture" the 
representation space $\Cal H(\sigma, \chi)$ consists of the closure of $C_c^\infty(V, W_\sigma)$ in 
the norm

$$\int\limits_V\Vert f(v)\Vert_\sigma^2e^{2\text{Re}\lambda(H(v)}dv$$
where $v=\kappa(v)$ exp $H(v)n(v)$ from the Iwasawa decomposition and $e^{2\text{Re}\lambda 
H(v)}=\vert\chi_\lambda
(\text{exp}~H(v))\vert^2.$  The $G$ action in this realization is given by

$$U(\sigma, \lambda)(g) 
f(x)=\sigma(m(g^{-1}x))^{-1}\chi_{\lambda+\rho}(a(g^{-1}x))^{-1}f(v(g^{-1}x))$$
where $g^{-1}x$ has been expressed in $VMAN$ coordinates.

Let $\chi_\lambda'$ be the contragredient character to $\chi_\lambda$ and form $\Cal H(\sigma, 
\lambda')$.  The mapping 
$f\to\chi_\lambda^2 (\text{exp}~H(x)) f(x)$ is an isometry of $\Cal H(\sigma, \lambda)$ onto $\Cal 
H(\sigma, \lambda')$.  
Moreover $\Cal H(\sigma, \lambda)$ and $\Cal H(\sigma, \lambda')$ are naturally dual via the pairing

$$(f, g)=\int\limits_V<f(v), g(v)>_\sigma dv.$$
In addition, for $f$ in $\Cal H(\sigma, \lambda), \, g~\text{in}~\Cal H(\sigma, \lambda')$ and 
$y~\text{in}~G,$ the 
matrix coefficient $c_{f,g}(y)$ is defined by

$$c_{f,g}(y)=(U(\sigma, \lambda)(y)f, g)$$
and one has the identity

$$(U(\sigma, \lambda)(y)f, g)=(f, U(\sigma, \lambda')(y^{-1})g).$$

The $A$-invariant measure on $V'$ is $\xi_{-\rho}(w^{-1}v) dv.$  Define a Hilbert space 
$\widetilde{\Cal H}(\sigma, 
\lambda)$ to consist of the closure of $C_c^\infty(V, W_\sigma)$ in the norm

$$\int\limits_V\Vert f(v)\Vert_\sigma^2e^{2\text{Re}\lambda H(v)} \xi_{-\rho}(w^{-1}v)dv.$$
The map $T:  \Cal H(\sigma, \lambda)\to\widetilde{\Cal H}(\sigma, \lambda)$ given by

$$Tf(x)=\xi_{-\rho}(w^{-1}x)^{-1/2}f(x), x\neq e$$
has inverse

$$T^{-1}f(x)=\xi_{-\rho}(w^{-1}x)^{1/2}f(x)$$
and establishes an isomorphism between these spaces.  There is then an equivalent realization of the 
principal series $\tilde{U}(\sigma, \lambda)(\cdot)$ on $\widetilde{\Cal H}(\sigma, \lambda)$ which 
we refer to as the ``$A$-adapted" realization and is given by

$$\tilde{U}(\sigma, \lambda)(\cdot)=T\circ U(\sigma, \lambda)(\cdot)\circ T^{-1}.$$
A straightforward computation using Lemma 3.1 and 3.2 gives

$$\lno{\tilde{U}(\sigma,\lambda)(a) f(x)=\chi_\lambda(a)\delta_a^\ast f(x), a\in A,&&(4.1)\cr}$$
$$\lno{\tilde{U}(\sigma, \lambda)(m) f(x)=\sigma(m)\delta_m^\ast f(x), m\in M.&&(4.2)\cr}$$
Furthermore, $\widetilde{\Cal H}(\sigma, \lambda)~\text{and}~\widetilde{\Cal H}(\sigma, \lambda')$ 
are in natural duality via
$$\lno{(f, g)=\int\limits_{V'}<f(x), g(x)>_\sigma\xi_{-\rho}(w^{-1}x) dx,&&(4.3)\cr}$$
and one has
$$\lno{(\tilde{U}(\sigma, \lambda)(y) f, g)=(f,\tilde{U}(\sigma, \lambda')(y^{-1})g).&&(4.4)\cr}$$

The $A$-adapted realization gives a new integral representation for matrix coefficients with close 
connections with abelian harmonic analysis. Using the construction from \S 2 we take $f_{\lambda}(t,\omega) = \chi_\lambda(t^{-1})f(t,\omega)$ with modified Fourier transform $\tilde {f}(\upsilon,\omega) =\hat{f}_{\lambda}(\upsilon,\omega)$, and similarly with the contragredient $\chi_\lambda'$ and $g$ we obtain $\tilde {g}$.

\proclaim{Theorem 4.1}  Let $\chi_\lambda$ be a unitary character of $A$ and let $\dim \frak 
a=\ell.$  Let $\sigma$ be a finite dimensional unitary representation of $M$. Let 
$f$ be in $\widetilde{\Cal H}(\sigma, \lambda)$ and $g$ be in $\widetilde{\Cal H}(\sigma, 
\lambda').$  Then, for $a$ in $A$

$$\lno{c_{f,g}(a)=\int\limits_{\Bbb R^\ell}e^{-i\upsilon\log a}m(\lambda, 
\upsilon)d\upsilon.&&(4.5)\cr}$$
Here $m(\lambda, \upsilon)$ is a function depending on $f$ and $g$ and given by

$$\lno{m(\lambda, \upsilon)=\int\limits_{A\bs V'}<\tilde{f}(\upsilon, \omega), 
\tilde{g}(\upsilon, \omega)>_\sigma d\omega.&&(4.6)\cr}$$\endproclaim

\remark{Remarks}  We shall refer to the measure $m(\lambda, \upsilon)d\upsilon$ as the Bochner 
measure 
for $c_{f,g}.$  Since $\tilde{U}$ is unitary the matrix coefficients $c_{f,f}$ are positive definite 
functions on $G$, a fortiori on $A$.  Hence Bochner's theorem assures the existence of a 
non-negative 
measure with this property while (4.5) shows it is absolutely continuous with respect to Haar 
measure and (4.6) gives a formula for it.\endremark

\demo{Proof}  We begin with (4.3)

$$\eqalignno{c_{f,g}(a)&=(\tilde{U}(\sigma, \lambda)(a) f,g)\cr
&=\int\limits_{V'}<\tilde{U}(\sigma, \lambda)(a) f(x), g(x)>_\sigma\xi_{-\rho}(w^{-1}x) dx\cr
&=\int\limits_{A\times A\bs V'}<\tilde{U}(\sigma, \lambda)(a) f(h,\omega), g(h,\omega)>_\sigma 
dh\times dw\cr
&=\int\limits_{A\times A\bs V'}<\chi_\lambda(a)\delta_a^\ast f(h,\omega), g(h,\omega)>_\sigma 
dh\times d
\omega\cr}$$
We continue with a formal application of Parseval's formula (2.3) and (2.2) to get

$$\eqalignno{c_{f,g}(a)&=\int\limits_{\hat{A}\times A\bs 
V'}\chi_\upsilon(a)^{-1}<\tilde{f}(\upsilon, \omega), 
\tilde{g}(\upsilon, \omega)>_\sigma d\upsilon\times d\omega\cr
&=\int\limits_{\Bbb R^\ell}e^{-i\upsilon\log a}m(\lambda, \upsilon)d\upsilon.\cr}$$
To justify the use of (2.3), we observe that $f$, $\tilde{U}f$ and $g$ have square integrable 
$\sigma$-norms with respect to $\xi_{-\rho}(w^{-1}x)dx$.  Hence, a.e. $d\omega$, $\tilde{U}f(\cdot, 
\omega)$ and $g(\cdot, \omega)$ have square integrable $\sigma$-norms relative to $dh$. Since $a\in A$ acts 
by $\chi_\lambda (a^{-1})$ and (2.3) is valid for each component of the vector the result follows. 
\hfill\bx \enddemo

\remark{Remarks}  The unitarity of $\chi_\lambda$ was used only to justify (2.3).  For non-unitary 
$\chi_\lambda$, whenever (2.3) can be used there is such an integral formula.  Also, notice that $f$ 
and $g$ 
need not be $K$-finite.\endremark
\vfil
\eject
\subhead \S 5.  Non-unitary principal series \endsubhead
\medskip

We assume in this section that dim $\frak a=1.$  Then $\frak a=\text{cl}(\frak a^+\cup w\frak 
a^+w^{-1}).$  
If $\lambda$ and $\upsilon$ are $\Bbb C$-valued linear forms on $\frak a$ we say $\lambda>\upsilon$ 
if Re$(\lambda
-\upsilon)(\frak a^+)>0.$

\proclaim{Theorem 5.1}  Let $\chi_\lambda$ be a character of $A$ with Re$\lambda\geq 0$ and both 
$\rho>\lambda$ and $\rho >-\lambda'.$ Let 
$\sigma$ 
be a finite dimensional unitary representation of $M$.  Let $f$ in $\widetilde{\Cal H}(\sigma, 
\lambda)$ and $g$ in 
$\widetilde{\Cal H}(\sigma, \lambda')$ be smooth vectors.  Then

$$\lno{c_{f,g}(a)=\int\limits_{\Bbb R} e^{-i\upsilon\log a}m(\lambda,\upsilon)d\upsilon, \,a\in 
A&&(5.1)\cr}$$
where

$$\lno{m(\lambda,\upsilon)=\int\limits_{A\bs V'}<\tilde{f}(\upsilon, \omega), \tilde{g}(\upsilon, 
\omega)>_\sigma 
d\omega&&(5.2)\cr}$$
\endproclaim

\demo {Proof}  As remarked in $\S 4$ it is enough to justify the use of (2.3).  Now it is well known 
that 
the intertwining operator between the compact and non-compact realization $\Cal H(\sigma, \lambda)$ 
is given by

$$\text{Ih}(x)=e^{-(\lambda+\rho)H(x)}h(\kappa(x)), \,x\in V.$$\enddemo
\noindent Following $I$ with $T$ we may suppose $f$ (and similarly $g$) is of the form

$$\xi_{-\rho}(w^{-1}x)^{-1/2}e^{-(\lambda+\rho)H(x)}h(\kappa(x))$$
with $h$ a $K$-finite $W_\sigma$-valued function on $K$.  In particular, each component of $h$ is 
bounded.  
Now, $\tilde{U}(\sigma, \lambda)(a) 
f(x)=\chi_\lambda(a)\xi_{-\rho}(w^{-1}a^{-1}xa)^{-1/2}e^{-(\lambda +\rho)
H(a^{-1}xa)}h(\kappa(a^{-1}xa)).$  Since $h$ has bounded components, for the purpose of estimates it 
may be 
ignored.  Choose a cross-section $A\bs V'\to V'$ and denote a fixed element in the range by 
$\omega$.  Recall 
the definition from $\S 2$ of $f_\lambda$ and $g_{\lambda'}$, i.e., $f_\lambda(t\omega 
t^{-1})=\chi_\lambda
(t)^{-1}f(t,\omega)$ and $g_{\lambda'}(t\omega t^{-1})=\chi_{\lambda'} (t)^{-1}g(t,\omega)\cdot$ We 
shall show that the Parseval formula (2.3) applies to the functions  
$\chi_\lambda(t)\xi_{-\rho}(w^{-1}t
\omega t^{-1})^{-1/2}e^{-(\lambda +\rho)H(t\omega t^{-1})}$ and 
 $\chi_{\lambda'}(t)\xi_{-\rho}(w^{-1}t\omega 
t^{-1})^{-1/2}e^{-(\lambda'+\rho)H(t\omega t^{-1})}$.

First notice that these functions are locally bounded in $t$ and so it suffices to consider their 
behavior for $t\to\infty$ in $A_+$ and $A_-=wA_+w^{-1}.$  Next, from Lemma 3.1 we have

$$\eqalignno{f_\lambda(t,\omega)
&\sim e^{-\lambda\log 
t}\chi_{-\rho}(w^{-1}tw)^{-1/2}\chi_{-\rho}(t^{-1})^{-1/2}\xi_{-\rho}(w^{-1}\omega)
^{-1/2}\cr
&\times e^{-(\lambda+\rho)H(t\omega t^{-1})}\cr}$$
and

$$\eqalignno{g_\lambda'(t,\omega)
&\sim e^{-\lambda'\log 
t}\chi_{-\rho}(w^{-1}tw)^{-1/2}\xi_{-\rho}(w^{-1}\omega)^{-1/2}e^{-(\lambda'+\rho)H
(t\omega t^{-1})}\cr}$$
and using Lemma 3.2 and ignoring $\xi_{-\rho}(w^{-1}\omega)^{-1/2}$ we get
$$f_\lambda(t,\omega)\sim e^{-\lambda\log t}e^{-\rho\log t}e^{-(\lambda+\rho)H(t\omega t^{-1})}$$
and

$$g_{\lambda'}(t,\omega)\sim e^{-\lambda'\log t}e^{-\rho\log t}e^{-(\lambda'+\rho)H(t\omega 
t^{-1})}.$$
Now, if $t\in A_+$, since $\omega$ is in 
$V',\lim\limits_{t\mathop{\rightarrow}\limits_{A_+}\infty}t\omega
t^{-1}=1.$  Thus for $t\mathop{\rightarrow}\limits_{A_+}\infty, \,e^{-(\lambda+\rho)H(t\omega 
t^{-1})}$ and $e^{-(\lambda'+\rho)H
(t\omega t^{-1})}$ are bounded, while $e^{-(\lambda+\rho)\log t}$ and $e^{-(\lambda'+\rho)\log t}$ 
are 
integrable on $A_+$ because $\rho >\lambda$ and also $\rho >-\lambda'.$

Before we consider $t$ in $A_{-}$ we shall express $f_\lambda$ and $g_{\lambda'}$ in another form.  
For 
$x$ in $V$

$$H(x)=H(w^{-1}x)=H(v(w^{-1}x))+\log a(w^{-1}x),$$
thus

$$H(t\omega t^{-1})=H(v(w^{-1}t\omega t^{-1})+\log a(w^{-1}t\omega t^{-1}).$$
Also, we have

$$\log a(w^{-1}t\omega t^{-1})=\log a(w^{-1}\omega)-\log t+\log(w^{-1}tw)$$
and

$$v(w^{-1}t\omega t^{-1})=w^{-1}tw ~v(w^{-1}\omega) w^{-1}t^{-1}w.$$
After discarding terms independent of $t$ we have

$$\eqalignno{f_\lambda(t,\omega)
&\sim e^{-(\lambda+\rho)\log t}e^{-(\lambda+\rho)\log(w^{-1}tw)}e^{(\lambda+\rho)\log t}\cr
\cr
&\times e^{-(\lambda+\rho)H(w^{-1}tw ~v(w^{-1}\omega)w^{-1}t^{-1}w)}\cr}$$
and

$$g_{\lambda'}(t,\omega)~\text{similarly~with}~\lambda\leftrightarrow\lambda'.$$
Simplifying gives

$$f_\lambda(t,\omega)\sim e^{-(\lambda+\rho)\log(w^{-1}tw)}e^{-(\lambda+\rho)H(w^{-1}tw 
~v(w^{-1}\omega)w^{-1}
t^{-1}w)}$$
and $g_{\lambda'}$ similarly with $\lambda\leftrightarrow\lambda'$.

Now let $t\in A_{-}.$  Then $w^{-1}tw\in A_+$ and Harish-Chandra [H-C] has shown that if $h\in A_+, 
\mu>0$ 
and $x\in V$ then

$$\mu(H(x)-H(hxh^{-1}))\leq 0.$$
So, with $\mu=\text{Re}~(\rho+\lambda)~\text{or}~\mu=\text{Re}~(\rho+\lambda')$ we get

$$\vert e^{-(\lambda+\rho)H(w^{-1}tw ~v(w^{-1}\omega)w^{-1}t^{-1}w)}\vert\leq 
e^{-\text{Re}(\lambda+\rho)H(w
^{-1}\omega)}$$
and similarly for the $\lambda'$ term.

Thus, ignoring terms bounded, for $t$ in $A_{-}$ we have 

$$\vert f_\lambda(t,\omega)\vert\leq e^{-\text{Re}(\lambda+\rho)\log(w^{-1}tw)}$$
and

$$\vert g_{\lambda'}(t,\omega)\vert\leq e^{-\text{Re}(\lambda'+\rho)\log(w^{-1}tw)}.$$
Both of which are integrable on $A_{-}.$  Hence, the Parseval relationship (1.3) can be applied to 
$f,g.$

\remark{Remark}  We note that the only use of dim $\frak a=1$ was to write $A=\text{cl}(A_+\cup 
A_{-}).$  
Otherwise, the estimates hold in general rank on $\frak a_+$ and the opposite chamber.\endremark

\subhead \S 6 Computation of Bochner measure 
\endsubhead
\medskip
The function part of the Bochner measure, $m(\lambda,\upsilon)$, contains 
much information about the matrix coefficients of unitary induced representations.  For complementary 
series the unitary structure on $\widetilde{\Cal H}(1,\lambda)$ must be modified by the intertwining 
operator, however on the $K$-fixed vector the intertwining operator in the ``compact picture" is a scalar. 
Thus we can use Theorem 5.1 for the spherical complementary series ($\sigma$ trivial) as well as for the 
spherical principal series  and the zonal spherical function.  To get some understanding 
of $m(\lambda,\upsilon)$ we shall compute it, exactly in these cases.

We first recall some calculations from [H2]. The set $\Sigma_+$ consists of, at most, $\{ \alpha, 
2\alpha\}.$  We let $p=\dim\frak g_\alpha$ and $q=\dim\frak g_{2\alpha}$ and choose the usual inner product 
on $\frak g$.  If $X$ denotes elements of 
$\frak g
_{-\alpha}$ and $Y$ those from $\frak g_{-2\alpha}$, each $x$ in $V$ can be written uniquely 
$x=\text{exp}~X~\text{exp}~Y.$  The relevant formulae 
are ([H2] p.59)

$$\lno{\xi_{-\rho}(w^{-1}x)=[c^2\vert X\vert^4+4c\vert Y\vert^2]^{-1/4(p+2q)}&&(6.1)\cr}$$

$$\lno{e^{\rho H(x)}=[(1+c\vert X\vert^2)^2+4c\vert Y\vert^2]^{1/4(p+2q)}&&(6.2)\cr}$$
where $c^{-1}=4(p+4q)$

As remarked in $\S 3, \frak v'=\frak v\bs \{0\}.$  It follows from Lemma 3.1 that the level set 
$\xi_{-\rho}
(w^{-1}x)=1$ provides a cross-section $A\bs V'\to V'.$  Also, being $\Bbb R$ rank 1 we express 
$\lambda=
\lambda\rho$ with $\vert\text{Re}~\lambda\vert<1$ in order to apply Theorem 5.1.

The $K$-fixed vector in $\widetilde{\Cal H}(1,\lambda)$ is 

$$\eqalignno{f_\lambda(a,\omega)&=e^{-\lambda\rho\log a}e^{-(\lambda\rho+\rho)H(a\omega 
a^{-1})}\xi_{-\rho}
(w^{-1}a\omega a^{-1})^{-1/2}\cr
\cr
&=e^{-\lambda\rho\log a}e^{(\lambda\rho+\rho)H(a\omega a^{-1})}e^{-\rho\log 
a}\xi_{-\rho}(w^{-1}\omega)
^{-1/2}.\cr}$$
Using $\xi_{-\rho}(w^{-1}\omega)=1$ and Helgason's formula (6.2) gives

$$f_\lambda(a,\omega)=e^{-(\lambda\rho +\rho)\log a}[(1+c\vert e^{-\alpha\log a}X\vert^2)^2+4c\vert
e^{-2\alpha\log a}Y|\vert^2]^{-{p+2q\over 2}(1+\lambda)}$$
\medskip
\noindent We set $s=e^{-\alpha\log a}$ and use $\rho=(p+2q)\alpha/2$ to get

$$f_\lambda(s,\omega)={s^{{p+2q\over 2}(1+\lambda)} \over [(1+cs^2\vert X\vert^2)^2+4cs^4\vert 
Y\vert^2]^
{{p+2q \over 2}(1+\lambda)}}.$$
\medskip
\noindent Next setting $t=s^2$ and using $1=\xi_{-\rho}(w^{-1}x)=[c^2\vert X\vert^4+4c\vert 
Y\vert^2]$ gives

$$f_\lambda(t,\omega)={{t^{{p+2q\over 4}(1+\lambda)}} \over {[1+2t~c\vert X\vert^2+t^2]}^
{{p+2q \over 4}(1+\lambda)}}.$$
\medskip
\noindent  Finally, since $c\vert X\vert^2\leq 1$ we write it as cos $\theta$.  Thus, using (2.1), 
we have

$$\eqalignno{\tilde{f}(\upsilon\rho,\omega)&=\int\limits_A\chi_{i\upsilon\rho}(a)^{-1}f_{\lambda}(a
,\omega)da\cr
&=\int\limits_0^\infty {t^{{p+2q\over 4}(1+\lambda-i\upsilon)} \over [1+2t\cos \theta +t^2]^{{p+2q 
\over 4}(1+\lambda)}}{dt\over t}\cr}$$ 
as the integral to be evaluated for $\lambda$ and similarly $\lambda'$.  Let $a=\big({p+2q\over 
4}\big)(1+\lambda-
i\upsilon)-1$ and $b=\big({p+2q\over 4}\big)(1+\lambda)$ and note that Re $a>-1$ and Re $b>0.$  We 
consider, under 
these circumstances, the integral

$$I=\int\limits_0^\infty{t^a\over [1+2t\cos \theta +t^2]^b}dt, \theta \in [0, {\pi\over 2}].$$
We factor the denominator and treat the integral as a contour integral.  Thus,

$$\eqalignno{I&=\int\limits_0^\infty{t^a \over (e^{i\theta}+t)^b(e^{-i\theta}+t)^b}dt\cr
\cr
&=e^{-2i~b\theta}e^{ia\theta}\int\limits_0^\infty{[e^{-i\theta} t]^a 
dt\over[1+e^{-i\theta}t]^b[e^{-2i\theta}+e^
{-i\theta}t]^b}.\cr}$$
 
After some elementary manipulations involving deformation of contours one arrives at

$$I=e^{i(a+1)\theta}(e^{2i\theta}-1)^{-b}\int\limits_0^1 x^a(1-x)^{2b-a-2}(x-z)^{-b} dx.$$
With $v=z^{-1}$ this is found in [W-W], p. 293 to be

$$I={\Gamma(a+1)\Gamma(2b-a-1) \over \Gamma(2b)}e^{i(a+1)\theta}F(b,a+1;2b;v),$$
using Gauss' hypergeometric function, provided $\vert v\vert<1.$  Now $\vert v\vert = \vert 
2\sin~\theta\vert<1$ means $0\leq\theta<{\pi\over 6}.$  To get the integral for $\theta\in[0, {\pi\over 
2}$] 
we 
analytically continue the hypergeometric function by the quadratic transformation ([BP])

$$F(\alpha, \beta; 2\beta; z)=(1-z)^{-\alpha\!/2}F({\alpha\over 2}, \beta-{\alpha\over 2}; 
\beta+{1\over 2}; {z^2 \over 
4(z-1)})$$
or, finally,

$$I={\Gamma(a+1)\Gamma(2b-a-1) \over \Gamma (2b)}F({a+1 \over 2}, b-{a+1 \over 2}; b+{1 \over 2}; 
\sin^2\theta).$$
Recalling how $a, b$ were defined we get

$$\lno{\tilde{f}(\upsilon, \omega)={\Gamma({(p+2q) \over 4} (1+\lambda-i\upsilon))\Gamma({(p+2q) 
\over 
4}(1+\lambda+i\upsilon)) \over \Gamma({(p+2q) \over 4}(2+2\lambda))}&&(6.3)\cr}$$

$$\times F(\alpha, \beta; \gamma;\sin^2\theta)~\text{where}$$

$$\eqalignno{\alpha&=\big({p+2q \over 4}\big)\big({1 \over 2}+{\lambda \over 2}-{i\upsilon \over 
2}\big), 
\beta=\big({p+2q \over 4}\big)\big({1 \over 2}+{\lambda \over 2}+{i\upsilon \over 2}\big)\cr
\cr
\gamma&=\big({p+2q \over 4}\big)\big(1+\lambda\big)+{1 \over 2}.\cr}$$
Simiarly, for $\lambda'$.

Some comments are in order.  When there is no $2\alpha$ root the $\cos \theta$ appears as 1 in the 
beginning and then the hypergeometric function is evaluated at zero, hence is 1 (a considerable 
simplification!).  Next, after we had evaluated the integral, we happened across it in a table.  
But, as similar computations for prehomogeneous vector spaces of parabolic type are a possibility, 
 we think it worth presenting.  Finally, compressing the notation, the resultant integral 
formula for the 
spherical function $\varphi_\lambda$,

$$\lno{\varphi_\lambda(a)=\int\limits_{\Bbb R}e^{-i\upsilon\log a}m(\lambda, 
\upsilon)d\upsilon&&(6.4)\cr}$$
seems new to us even from the special function literature.  The only analogous, though definitely 
different, integrals were obtained by Barnes [Ba].
\vfil
\eject

\subhead \S 7.  Applications \endsubhead

In this section, we give four applications of the previous sections.  

\subsubhead Positive Definite \endsubsubhead

We have already remarked that $m(\lambda,\upsilon)\geq 0$ 
is necessary for unitarity of $\tilde{U}.$  Conversely, one might ask if it is sufficient, or, in 
other words, if the restriction to $A$ are positive definite must $\tilde{U}$ be unitary?

\proclaim{Proposition 7.1}  Let $\lambda=\lambda\rho$ and suppose $\vert\text{Re}~ \lambda\vert\leq 
1$.  If 
either (i) Re $\lambda=0$ or (ii) Im $\lambda=0,$ then $\varphi_\lambda$ is positive definite on A.
\endproclaim

\demo {Proof}  First suppose $\vert\text{Re}~\lambda\vert<1$ for then we may use (6.3).  Of course 
it 
suffices to show that $m(\lambda,\upsilon\geq 0).$  If Re $\lambda=0$, then $\lambda'=\lambda$ and 

$$m(\lambda,\upsilon)=\int\limits_{A\bs V'}\tilde{f}(\upsilon,\omega)\overline{\tilde{f}(\upsilon, 
\omega)} 
dw\geq 0.$$
If Im $\lambda=0,$ then $\lambda'=-\lambda.$  In (6.3) the parameters $\alpha, \beta$ in the 
hypergeometric function are complex conjugates, while $\gamma$ is positive, so the function is 
positive.  
Also, the $\Gamma$-functions in the numerator are conjugates while the denominator is positive.  
Hence, 
$\tilde{f}(\upsilon, \omega)>0$ for $\lambda$ and similarly for $\lambda'.$  If $\lambda=\rho$, a 
careful look at the integral I shows the calculation is valid for $\lambda$ and as just observed is 
a 
positive function.  But for $\lambda'=-\rho$ we must compute the Fourier transform of 1, getting a 
Dirac measure.  Thus, the Bochner measure for $\lambda=\rho$ is  positive multiple of a Dirac 
measure at 
$\upsilon=0.$ \enddemo

In [K] the unitarizable spherical principal series are determined and shown to be (i) Re $\lambda=0$ 
or 
(ii) Im $\lambda=0$ and $\vert \lambda\vert\leq {p\over 2}$ when $q=0$ and when $q\ne 0$ one also has 
${{p+2q}\over 2}$.  Comparison with Proposition 7.1 suggests an extension problem. 

\proclaim{Problem 1}  Let $f$ be defined on $G$ and suppose $f$ is positive definite on $A$.  Find 
sufficient conditions for $f$ to be positive definite on $G$.
\endproclaim

If instead we restricted to a larger group there is an easy answer.

\proclaim {Lemma 7.2}  Let $G$ be a locally compact group, $K$ a closed subgroup $S\subseteq G$ with 
$G=KS.$  Let $f$ be defined on $G$ and bi-$K$-invariant.  Then, $f$ is positive definite on $G$ if 
and 
only if $f$ is positive definite on $S$.
\endproclaim

\demo {Proof}  It is enough to show $f$ is p.d. on $G$.  Let $x_1,\cdots, x_n$ be in $G$ and say 
$x_i=k_is_i$  Let $c_i,\cdots,c_n$ be complex numbers.

$$\lno{\sum_{i.j}f(x_ix_j^{-1})\overline{c}_ic_j&=\sum_{i,j} f(k_is_is_j^{-1}k_j^{-1})\overline{c_i} 
cj\cr
&=\sum_{i,j}f(x_is_j^{-1})\overline{c_i}c_j\geq 0.\cr}$$\enddemo

In [F-K] the set $S=AK$ was used and the $K$-type expansion together with Bochner's theorem on the 
compact group $K$ gave another proof of Kostant's theorem.  Alternatively $S=MAN$ leads to the use 
of 
Kirillov theory on $N$ together with a Bochner theorem for nilpotent groups.

\subsubhead Harish-Chandra expansion
\endsubsubhead

We shall obtain from the integral formula (6.4) another derivation of Harish-Chandra's expansion for the 
elementary spherical function $\varphi
_\lambda(a)$ but without the use of any differential equations. Fix $\text {H}\in \frak a$ of norm one, and 
set $a_t = exp tH$. Recall that the level set $\Omega = \{\xi_{-\rho}(w^{-1}x = 1\}$ is a cross-section for 
$A\bs V'$. From (5.1) 
$$\varphi_\lambda (a_t) =\int\limits_{\Bbb R} e^{-i\upsilon t} \int\limits_{\Omega}<\tilde{f}(\upsilon, 
\omega), \tilde{g}(\upsilon, \omega)>_\sigma d\omega d\upsilon $$
here $\tilde{f}(\upsilon, \omega)$ is given by (6.3). 
But first, we need 

\proclaim {Lemma 7.3}  If $\lambda=\lambda\rho$ with $\vert\text{Re}~\lambda\vert<1$ then the 
Bochner 
measure is of the form

$$m(\lambda, \upsilon)=\Upsilon(\lambda, \upsilon)\Upsilon(-\lambda, \upsilon)h(\lambda,\upsilon)$$
where $h$ is entire in $\upsilon\in\frak a_c$ and

$$\Upsilon(\lambda, \upsilon)={\Gamma(r(1+\lambda-i\upsilon))\Gamma(r(1+\lambda+i\upsilon)) \over 
\Gamma(r(2+2\lambda))}$$
here $r=(p+2q)/4$. \endproclaim

\demo{Proof}  From (6.3) we define $\Upsilon(\lambda,\upsilon)$ to be the gamma functions appearing 
there.  
Next, simply observe that the complex conjugate of (4.5) for $\lambda'$ simply takes 
$\lambda\to-\lambda$.  
We must consider the integral over $A\bs V'$ of the product of two hypergeometric functions.  But, 
the 
$\upsilon$-dependence is only in the $\alpha$ and $\beta$ parameters in which the hypergeometric 
function 
is entire (note that Re$(\Upsilon-\alpha-\beta)>0)$ and as above, the complex conjugate can be 
absorbed as 
a change of sign since $A\bs V'$ may be viewed as the compact set $\xi_{-\rho}(w^{-1}x)=1$ we have 
the 
resulting integral, $h(\lambda,\upsilon)$, entire in $\upsilon \in \frak a_c$.\enddemo

We illustrate the computation by giving the details for $\text {Sl}(2, \Bbb R)$ and $\lambda$ imaginary. In order to compare the result to ones in the literature we change notation slightly. Henceforth $\lambda$ imaginary will be written as $i\lambda :=i\lambda \alpha$, and similarly in the integral we replace $\upsilon$ with $\upsilon \alpha$ and $tH$ with the co-root $t\check {H_\alpha}$. In 
this case the integral becomes with $r=1/4$.
$$\lno{&&(7.1)\cr} $$
$$\varphi_\lambda (a_t) =c\int\limits_{\Bbb R} e^{-i\upsilon t} 
{\Gamma(r(1+i({\lambda\over 2}-{\upsilon\over 2})))\Gamma(r(1+i({\lambda\over 2}
+{\upsilon\over 2})))\over 
\Gamma(r(2+4\lambda))} {\Gamma(r(1-i({\lambda\over 2}+{\upsilon\over 2})))\Gamma(r(1-i({\lambda\over 2}-{\upsilon\over 2})))\over 
\Gamma(r(2-4\lambda))}d\upsilon$$

\noindent The integral consists of terms of the form $\Gamma(r \pm i(r/2)(\lambda\pm\upsilon))$ and these 
functions have poles at $\pm \lambda\pm i(2k+1/2)$. For $\lambda \ne 0$ all poles of the integrand are 
simple and one easily sees that the residues are of $\Gamma (z)$ when $z = -k$. To compute $\text 
{Res}_{-k} \Gamma (z)$ use the duplication formula $\Gamma (z)\Gamma (1-z) = {\pi \over \text{sin}\pi z}$ 
and one obtains 
$$\text {Res}_{-k} \Gamma (z) = {(-1)^k\pi \over \Gamma(k+1)}.$$

\noindent To evaluate the other terms one needs to simplify $\Gamma(-k+i\lambda)$. Again, a use of the 
duplication formula and the fact that $\text {sin }\pi i\lambda = {(-1)\pi\over {i\lambda 
\Gamma(i\lambda)\Gamma(-i\lambda)}}$ gives 

$$\Gamma(-k+i\lambda) ={(-1)^{k+1}i\lambda \Gamma(i\lambda)\Gamma(-i\lambda)\over \Gamma(1+k-i\lambda)}.$$

\noindent We evaluate (7.1) by an application of Cauchy's theorem and a rectilinear contour from -R to R to 
 R+iR to -R+iR to -R. It is routine to show the contributions other than along -R to R tend to zero as 
$R\to\infty$. The result is

$$\eqalignno{ \varphi_\lambda (a_t) &= {\Gamma(i\lambda)\over\Gamma(1/2 
+i\lambda)}\sum_{k=0}^{\infty}e^{-i\lambda t}e^{-(2k+1/2)t}{\Gamma(k+1/2)\over \Gamma(k+1)}{\Gamma(1/2 
+k-i\lambda)\over \Gamma(1/2 -i\lambda)}{\Gamma(1-i\lambda)\over\Gamma(1+k+i\lambda)}\cr
 &+ {\Gamma(-i\lambda)\over\Gamma(1/2 -i\lambda)}\sum_{k=0}^{\infty}e^{i\lambda 
t}e^{-(2k+1/2)t}{\Gamma(k+1/2)\over \Gamma(k+1)}{\Gamma(1/2 +k+i\lambda)\over \Gamma(1/2 
+i\lambda)}{\Gamma(1+i\lambda)\over\Gamma(1+k-i\lambda)}.\cr}$$

\remark{Remarks} (1) It is curious that this derivation uses in no explicit way the differential equation 
for the spherical function, in contrast to all other derivations of it.

(2) Obviously the same computation works for $\text {SO }(n,1)$ using a different value of ``$r$". The computation 
for the other rank one groups is complicated by having to evaluate 

$$\int\limits_0^{\pi/2}\vert F(\alpha,\beta;\gamma; sin^2(\theta))\vert^2 d\theta.$$

Our only progress to obtain a clean formula here was to expand one of the hypergeometric functions and use 
[B] p.399 to evaluate each of the resulting integrals; thereby obtaining an infinite series in quotients of 
$\Gamma$ functions at various arguments. This was sufficient to obtain the "c" functions and the leading 
term.
\endremark

\subsubhead Holomorphic extension
\endsubsubhead

Again, we illustrate the method by presenting the details solely for $\text{Sl}(2,\Bbb R)$. We begin with (7.1).
$$\varphi_\lambda (a_t) =c\int\limits_{\Bbb R} e^{-i\upsilon t} 
{\Gamma(r(1+i({\lambda\over 2}-{\upsilon\over 2})))\Gamma(r(1+i({\lambda\over 2}
+{\upsilon\over 2})))\over 
\Gamma(r(2+4\lambda))} {\Gamma(r(1-i({\lambda\over 2}+{\upsilon\over 2})))\Gamma(r(1-i({\lambda\over 2}-{\upsilon\over 2})))\over 
\Gamma(r(2-4\lambda))}d\upsilon.$$
We are interested in the Paley-Wiener phenomenon, i.e. exponential decay of the Fourier transform is 
equivalent to a holomorphic extension of the function to a tube. The case $\lambda = 0$ contains all the 
essential technicalities; so we present that. Hence consider
$$I =  \int\limits_{\Bbb R} e^{i\upsilon x} \vert\Gamma({1\over 4}- {i\over 2}\upsilon)\vert^2 
\vert\Gamma({1\over 4}+{i\over 2}\upsilon)\vert^2 d\upsilon.$$

\noindent Since $\ov{\Gamma(u+iv)} =\Gamma(u-iv)$ the $\Gamma$ functions in the integral can be rewritten 
as 
$$[ \Gamma({1\over 4}+{i\over 2}\upsilon)\Gamma({1\over 4}-{i\over 2}\upsilon)]^2.$$
Take the logarithm and use Binet-Stirling to get
$$\eqalignno{\text{log }\Gamma({1\over 4}+{i\over 2}\upsilon)
 + \text{log }\Gamma({1\over 4}-{i\over 2}\upsilon) &= -{1\over 2} \text{log }\vert {1\over 4} + {i\over 
2}\upsilon\vert -{\upsilon\over 2}\text{arg }({1\over 4}+{i\over 2}\upsilon)\cr
 & +{\upsilon\over 2}\text{arg }({1\over 4}-{i\over 2}\upsilon)-{i\over 4}\text{arg }({1\over 4}+{i\over 
2}\upsilon)\cr
 &-{i\over 4}\text{arg }({1\over 4}-{i\over 2}\upsilon) +O(1).\cr}$$
\noindent Now for $\upsilon>0$ we have $\text{arg }({1\over 4}\pm {i\over 2}\upsilon)\sim \pm 
\pi/2\mp{1\over 2\upsilon}, \upsilon\to\infty$, while for $\upsilon<0$ we have $\text{arg }({1\over 4}\pm 
{i\over 2}\upsilon)\sim \mp \pi/2\pm{1\over 2\vert\upsilon\vert}, \upsilon\to\infty.$ Then 

$$\text{arg }({1\over 4}+{i\over 2}\upsilon) - \text{arg }({1\over 4} -{i\over 2}\upsilon)\sim 
\pm\pi\mp{1\over\vert\upsilon\vert },\, \upsilon<0 (\,resp. \,\upsilon>0),$$
while 

$$\text{arg }({1\over 4}+{i\over 2}\upsilon) + \text{arg }({1\over 4} -{i\over 2}\upsilon) = o(1), \upsilon 
\to\infty.$$

\noindent Then
$$\eqalignno{\text{log }\Gamma({1\over 4}+{i\over 2}\upsilon) + \text{log }\Gamma({1\over 4}-{i\over 
2}\upsilon) &\sim \text{log } \vert {1\over 4} + {i\over 2}\upsilon\vert^{-1/2} -{\upsilon\over 2} [ 
\pm\pi\mp{1\over\vert\upsilon\vert}]+O(1)\cr
&\sim\text{log }\vert {1\over 4} + {i\over 2}\upsilon\vert^{-1/2} -\vert\upsilon\vert\pi/2 +O(1).\cr}$$
\noindent Consequently 
$$[ \Gamma({1\over 4}+{i\over 2}\upsilon)\Gamma({1\over 4}-{i\over 2}\upsilon)]^2 \sim 
c{e^{-\pi\vert\upsilon\vert}\over\sqrt{{1\over 16} +{\upsilon^2\over4}}}.$$
So the spherical function holomorphically extends in ``$x$" to ``$x+iy$". Moreover, one can detect the 
nature of the singularity. Indeed, taking $x=0$ and writing $y=\pi -\epsilon$, we obtain the Laplace 
transform
$$\int\limits_{\vert\upsilon\vert>R} e^{-\epsilon\upsilon}{d\upsilon\over\sqrt{{1\over 16} 
+{\upsilon^{2}\over 4}}}$$
which is easily seen to be $\sim \text{ln }\epsilon$. Both these results can be found in [KS-I] with a 
different proof. In fact, the preceeding computation was used as independent verification of the 
computations therein.
\remark{Remarks} The extension to $\lambda \ne 0$ presents no difficulties, nor does the extension to 
$SO(n,1)$. The remaining rank 1 groups can be done this way; albeit with additional work. As usual, the 
sticky point is (5.2) which here is the integral over $\Omega$ of $ F(2a, 2b; a+b+{1 \over 2}; 
\sin^2{\theta\over2})$. However we need only estimate - not evaluate - it. After some transformations on 
the parameters of the hypergeometric function one obtains (up to easily estimated factors) the Legendre 
function $P_n^m(cos\phi)$ where $n\in \Bbb C, \,\vert n\vert\to\infty$. For this the asymptotics in [Wa] 
p.291 are adequate to obtain estimates which allow one to detect the size of the  holomorphic extension, 
and the nature of the singularity. We omit the details because there is a detailed proof of the result in 
[KS-I].
\endremark

\subsubhead Spectral $F_f$
\endsubsubhead

We begin by recalling properties of the Abel transform as defined by Harish-Chandra. For $f\in C_c(G//K)$ 
define $F_f$ by 
$$F_f(a) = e^{\rho log a }\int\limits_N f(an) dn.$$
The Abel transform $\Cal A:f\to F_f$ has several interesting properties of which we highlight only those 
immediately relevant.
$$\lno{\Cal A: C_c(G//K)\longrightarrow C_c^W(A)&&(7.2)\cr}$$
where $C_c^W(A)$ denotes the compactly supported Weyl group invariant functions on $A$.
$$\lno{\Cal A \text{ is continuous in the usual topologies.}&&(7.3)\cr}$$
(If one uses the Harish-Chandra-Schwarz space, then an analogous version of both (7.2) and (7.3) are valid.)
Hence there exists $\Cal A': C_c^W(A)'\longrightarrow C_c(G//K)'$ satisfying for $f\in C_c(G//K)$ and $T\in 
C_c^W(A)'$
$$<\Cal Af,T> = <f,\Cal A'T>.$$
Let $\varphi_\lambda$, $\lambda$ pure imaginary, be the elementary spherical function. Then we have that 
$\varphi_\lambda\in C_c(G//K)'$; similarly we have the averaged character ${1\over\vert 
W\vert}\sum\limits_W e^{w\lambda(\cdot)} \in C_c^W(A)'$. Using $\tilde{f}$ to denote the spherical 
transform of $f$ we recall the familiar identity 
$$\lno{\tilde{f}(\lambda) = \hat {F_f}(\lambda)&&(7.4)\cr}$$

The duality relationship 
$$<\Cal Af,{1\over\vert W\vert}\sum\limits_W e^{w\lambda(\cdot)}> = <f, \Cal A'({1\over\vert 
W\vert}\sum\limits_W e^{w\lambda(\cdot)})>,$$
combined with the fact that $ \Cal A' \text{  is 1-1 }$ give the identity  
$$\lno{\varphi_\lambda (\cdot) = \Cal A' {1\over\vert W\vert}\sum\limits_W e^{w\lambda(\cdot)}.&&(7.5)\cr}$$
Indeed this relationship can be read in reverse, so that (7.4) and (7.5) are equivalent. As an aside, in [S-T] the specific form, in the case of real rank 1, of $\Cal A'$ as an integral operator in (7.5) was used to derive an expansion for the elementary spherical functions in terms of spherical functions for the associated Cartan motion groups.

Now take $f\in \Cal S(G//K)$ or just in $C_c(G//K)$ but with sufficient decay of its spherical Fourier transform, $\tilde{f}$. Recall the Bochner measure $m(\lambda, \upsilon)$. Define $F^s_{\tilde{f}}(\upsilon)$ by
$$\lno{ F^s_{\tilde{f}}(\upsilon) = \int\limits_{i\frak a^*}\tilde{f}(\lambda)m(\lambda, \upsilon) {d\lambda\over \vert c(\lambda)\vert^2}.&&(7.6)\cr}$$
Recall the inversion formula for the spherical transform
$$f(a) = \int\limits_{i\frak a^*}\tilde{f}(\lambda)\varphi_\lambda(a){d\lambda\over \vert c(\lambda)\vert^2}.$$
Then an elementary application of the Fubini theorem and the Euclidean Fourier inversion theorem gives
$$\lno{f(a) = {F^s_{\tilde{f}}}^{\vee}(a)&&(7.7)\cr}$$
which is to be contrasted with (7.4).

\vfil
\eject

\enddocument